\documentclass[12pt]{amsart}

\usepackage{amssymb,verbatim}

\theoremstyle{plain}
\newtheorem{thm}{Theorem}
\newtheorem{prop}{Proposition}
\newtheorem{cor}{Corollary}
\newtheorem{lem}{Lemma}

\theoremstyle{definition}
\newtheorem{de}{Definition}
\newtheorem{ex}{Example}
\newtheorem{rk}{Remark}

\renewcommand{\iff}{if and only if}

\newcommand{\hwc}{horizontally weakly conformal}
\newcommand{\h}{\mathcal{H}}
\renewcommand{\v}{\mathcal{V}}
\newcommand{\rn}{{\mathbb R}}
\newcommand{\sn}{{\mathbb S}}
\newcommand{\hn}{{\mathbb H}}

\newcommand{\bn}{{\mathbb B}}
\newcommand{\co}{\circ}

\DeclareMathOperator{\tr}{trace} \DeclareMathOperator{\g}{grad}
\DeclareMathOperator{\di}{div} 
\DeclareMathOperator{\ricci}{Ricci}

\begin{document}

\begin{abstract}
Inspired by the all-important conformal invariance of harmonic maps on
two-dimensional domains, this article studies the relationship between biharmonicity
and conformality.  We first give a characterization of biharmonic morphisms,
analogues of harmonic morphisms investigated by Fuglede and Ishihara, which, in
particular, explicits the conditions required for a conformal map in dimension four
to preserve biharmonicity and helps producing the first example of a biharmonic
morphism which is not a special type of harmonic morphism. Then, we compute
 the bitension field of horizontally weakly conformal maps, which
include conformal mappings. This leads to several examples of proper (i.e. non-harmonic) biharmonic
 conformal maps, in which dimension four plays a pivotal role.
We also construct a family of Riemannian submersions which are proper biharmonic maps.
\end{abstract}

\title[Biharmonic maps and morphisms]{Biharmonic maps and morphisms from conformal mappings}
\author{E. Loubeau}
\address{D{\'e}partement de Math{\'e}matiques\\
Universit{\'e} de Bretagne Occidentale\\
6, avenue Victor Le Gorgeu\\
BP 809\\
29285 Brest Cedex\\
France.} \email{Eric.Loubeau@univ-brest.fr}
\author{Y.-L. Ou*}
\address{Department of Mathematics \\
Texas A\&M University-Commerce \\
Commerce, TX 75429\\
USA.} \email{yelin\_ou@tamu-commerce.edu}
\thanks{Part of this work was done at the University of California at Riverside. \\
* Supported by the Texas A \& M University-Commerce Faculty Mini-Grant Program 2007.}
\date{\today}
 \maketitle

\section{Introduction}

A central feature of harmonic maps is their conformal invariance in dimension two.
Not only this allows defining harmonic maps on Riemann
surfaces but it also is the starting point to many properties of minimal branched immersions.

In the higher-order theory of biharmonic maps, one could expect similar properties in
dimension four. While this dimension certainly enjoys a special role for
biharmonicity, as illustrated by the conformal deformation of harmonic maps
of~\cite{B-K}, the study of the biharmonic stress-energy tensor or the
characterization of biharmonic morphisms of Section~\ref{section1}, no conformal
invariance of any sort has ever been observed.

Nevertheless, the interaction between conformality and biharmonicity remains a rich
subject and provides an interesting source of new examples.

In the theory of harmonic maps, a particularly fruitful approach has been to consider
maps which preserve local harmonic functions, called harmonic morphisms, because
their characterization as horizontally weakly conformal (a generalization of
Riemannian submersions) harmonic maps confers them a more geometrical flavour, which
counterweighs the analytical nature of harmonic maps (\cite{Fug,Ish}).

Their numerous geometrical properties have earned harmonic morphisms a choice place
among harmonic maps. The counterparts of these are biharmonic morphisms (see~\cite{OuB}),
 maps which pull back local biharmonic functions onto biharmonic
functions (and also, as it turns out, maps) and we characterize them as horizontally
weakly conformal biharmonic maps which are $4$-harmonic and satisfy an additional
equation, whose significance remains largely enigmatic.

While the number of conditions is directly due to the order of the problem, the
appearance of $4$-harmonicity is yet another clue to the specificity of dimension
four.

This characterization also solves one aspect of the original question on the
conformal invariance of biharmonic maps on four dimensional domains, since two extra
conditions are required.

\section{$p$-Harmonicity and biharmonicity}

At the origin of this work, lie the search and study of maps selected as extremals of
a measured quantity and the most natural class of functionals on the space of maps
between Riemannian manifolds is the $p$-energies.
\begin{de}
Let $\phi : (M,g) \to (N,h)$ be a smooth map between Riemannian manifolds and assume
$(M,g)$ compact then, for $p\in \rn$ ($p>1$), its {\em $p$-energy} is
$$ E_{p}(\phi) = \tfrac{1}{p} \int_{M} |d\phi|^p \, v_{g} .$$
The critical points of $E_{p}$ are called {\em $p$-harmonic maps} and simply {\em harmonic maps} for $p=2$.\\
Standard arguments yield the associated Euler-Lagrange equation, the vanishing of the
{\em $p$-tension field}:
$$ \tau_{p}(\phi) = |d\phi|^{p-4} [|d\phi|^{2}\tau(\phi) + \tfrac{p-2}{2} d\phi(\g |d\phi|^2)] =0  ,$$
where
$$ \tau(\phi) =  \tr{ \nabla d \phi } $$
is the tension field.
\end{de}

\begin{rk}
i) A fundamental feature of $p$-harmonicity for $p\neq 2$, is the collapse of
ellipticity at critical points and its negative consequences for regularity
properties (cf.~\cite{H-L}). For $p=2$, smoothness of continuous
harmonic maps is ensured by boot-strap methods and strong uniqueness, for example, follows.\\
ii) The main existence result for $p$-harmonic maps ($p>2$) is due to Duzaar and
Fuchs in~\cite{D-F} and generalizes the Eells-Sampson theorem.\\
iii) Working directly with the functional, one sees that, if $p=\dim{M}$, $E_{p}$ is
conformally invariant, a situation that seems more natural since, when the target
has trivial $p$-th homotopy group,
 $p$-harmonic maps exist in each homotopy class (\cite{Jost}).
\end{rk}

While harmonic maps have been extensively studied and shown to exist in numerous
circumstances (cf.~\cite{E-L}), in some situations they cannot exist or are very
limited. It is therefore interesting to turn to an alternative measuring the default
of harmonicity.
\begin{de}
Let $\phi : (M,g) \to (N,h)$ be a smooth map between compact Riemannian manifolds.
Define its {\em bienergy} as
$$ E^{2}(\phi) = \tfrac{1}{2} \int_{M} | \tau(\phi) |^{2} \, v_{g} .$$
Critical points of the functional $E^{2}$ are called {\em biharmonic maps} and its
associated Euler-Lagrange equation is the vanishing of the {\em bitension field}
$$ \tau^{2}(\phi) = - \Delta^{\phi} \tau(\phi) - \tr_{g}{R^{N}}(
d\phi, \tau(\phi) ) d\phi ,$$ where $\Delta^{\phi} = - \tr_{g} (\nabla^{\phi}\nabla^{\phi} - \nabla^{\phi}_{\nabla})$ is the
Laplacian on sections of $\phi^{-1}TN$ and $R^{N}$ the Riemann curvature operator on
$(N,h)$.
\end{de}

\begin{rk}
i) Clearly harmonic maps are automatically biharmonic, actually absolute minimums of
$E^2$.
For compact domains and negatively curved targets, the converse holds (\cite{Jiang}).\\
ii) An alternative to $E^2$ is to view $(N,h)$ isometrically immersed in $\rn^N$ and,
considering $\phi$ as a vector, take the $L^2$-norm of $\Delta\phi$ and call the
critical points {\em (extrinsic) biharmonic maps}. Except for flat targets, the two
definitions are distinct.
The regularity of both types of biharmonic maps has been extensively studied in~\cite{C-W-Y,W1,W2,W3}.\\
iii) If $M$ is non compact, we extend all these definitions by integrating over
compact subsets.
\end{rk}

A natural generalization of isometric immersions are weakly conformal maps, i.e.
whose differential either vanishes or is injective and conformal. Of particular
interest is their harmonicity since it characterizes minimality of the image and
defines minimal branched immersions. They also explain the conformal invariance of
harmonic maps on surfaces, directly at the level of the tension field.

Dual to this notion is horizontal weak conformality where, pointwise, the
differential is required to vanish or be surjective and conformal. This reverses the
constraint on the dimensions and enables a preservation of harmonicity in higher
dimensions.

\begin{de}
Let $\phi : (M^m ,g) \to (N^n ,h)$ be a smooth map between Riemannian manifolds. For
any point $x\in M$, let $\v_{x} = \ker d\phi_{x}$ be the vertical space at $x$ and
$\h_{x} = (\v_{x})^{\perp}$ the horizontal space at $x$. These spaces define a
vertical and a horizontal distribution. The map $\phi$ is called {\em horizontal
weakly conformal} if for any point $x\in M$ either $d\phi_{x} =0$ or $d\phi_{x}$ is surjective and
conformal from $\h_{x}$ to $T_{\phi(x)}N$, i.e.
$$ h( d\phi_{x}(X), d\phi_{x}(Y)) = \lambda^{2}(x) g(X,Y), $$
for all $X,Y \in \h_{x}$ and the function $\lambda$ is called the dilation of $\phi$.\\
For a vector $X\in TM$, $X^\h$ and $X^\v$ will denote the horizontal and vertical
parts of $X$. Points where $d\phi_{x} \neq 0$ are called {\it regular}.
\end{de}

\begin{rk}
i) If $m<n$, then $\phi$ is horizontally weakly conformal if and only if it is constant. \\
ii) By extending the function $\lambda$ by zero over critical points, we obtain a
smooth
function $\lambda^2$ defined on the whole of $M$. Besides $|d\phi|^2 = n \lambda^2$.\\
iii) Harmonic morphisms, i.e. maps which preserve local harmonic functions by
composition on the right-hand side, were characterized by Fuglede and Ishihara, as
horizontally weakly conformal harmonic maps (\cite{Fug,Ish}).
\end{rk}

\noindent Conventions: We will systematically use the Einstein convention on summing repeated indices.\\
Our convention for the Riemann curvature tensor will be
$$ R(X,Y) = [ \nabla_{X} , \nabla_{Y} ] - \nabla_{[X,Y]} ,$$
and the Laplacian on functions has been chosen with negative eigenvalues, i.e. $
\Delta f = \tr \nabla df$ but on vector fields $\Delta X = - \tr \nabla^{2} X$.

\section{The characterization of biharmonic morphisms}\label{section1}

In this section, we give an improvement of the characterization of biharmonic
morphisms and show that the inversion in the unit sphere $\rn^n \setminus \{0\} \to
\rn^n$ is a biharmonic morphism if and only if $n=4$, thus providing an example
which, unlike any of the previously known ones, is not a harmonic morphism.

In light of the theory of harmonic morphisms, cf. the monograph~\cite{B-W}, it only
seemed natural to define and study their biharmonic counterparts.

\begin{de}
A continuous map $\phi : (M,g) \to (N,h)$ is a {\em biharmonic morphism} if, for any
biharmonic function $f : U\subset N \to \rn$, such that $\phi^{-1}(U)\neq \emptyset$,
the pull-back function $f\circ\phi : \phi^{-1}(U) \subset M \to \rn$ is also
biharmonic.
\end{de}

\begin{rk}
i) Clearly constant maps and isometries are biharmonic morphisms and the composition
of
two biharmonic morphisms is again a biharmonic morphism.\\
ii) Since harmonic functions are automatically biharmonic, a biharmonic morphism will
pull-back harmonic functions onto biharmonic ones, but not necessarily harmonic.
Likewise, there is no reason to believe that a harmonic morphism should be a
biharmonic morphism. This distinction is clarified by Theorem~\ref{thm1} (see
also~\cite{L-O}).
Nevertheless, all the previously known examples of biharmonic morphisms came from special types of harmonic morphisms (\cite{L-O,OuB,Ou}).\\
iii) The existence of local harmonic coordinates on Riemannian manifolds
(\cite{D-K}), implies that in such a coordinate system, the components of a
biharmonic morphism are continuous biharmonic functions, hence smooth by standard
properties of elliptic partial differential equations. Therefore a biharmonic
morphism is always smooth.
\end{rk}

The geometric method used by Ishihara in~\cite{Ish} to characterize harmonic
morphisms as horizontally weakly conformal harmonic maps, and co-opted for the
semi-Riemannian case in~\cite{Fug96} by Fuglede (who had his own approach for
Riemannian metrics) is extended here to biharmonicity. As one should expect, four
conditions are needed to describe biharmonic morphisms, with the last one still very
much unfathomable.

\begin{thm} \label{thm1}
Let $\phi : (M^{m},g) \to (N^{n},h)$ be a smooth map between Riemannian manifolds.
Then $\phi$ is a biharmonic morphism {\iff} it is a {\hwc} biharmonic $4$-harmonic
map, of dilation $\lambda$, such that
\begin{align}\label{longeq}
& |\tau(\phi)|^{4} - 2\Delta \lambda^{2} |\tau(\phi)|^{2} + 4 \Delta \lambda^{2}
\di{\langle d\phi , \tau(\phi)\rangle} + n (\Delta \lambda^{2})^{2} \\
& + 2 \langle d\phi , \tau(\phi) \rangle (\nabla |\tau(\phi)|^{2}) + |S|^{2} = 0 ,\notag
\end{align}
where $S \in \odot^{2} \phi^{-1}TN$ is the symmetrization of the $g$-trace of $d\phi
\otimes \nabla^{\phi} \tau(\phi)$ and $\langle d\phi , \tau(\phi)\rangle (X) = \langle d\phi(X) , \tau(\phi)\rangle$.
\end{thm}

\begin{proof}
Let $\phi : (M^{m},g) \to (N^{n},h)$ be a smooth map between Riemannian manifolds,
then the statement of Theorem~\ref{thm1} is that $\phi$ is a biharmonic morphism
{\iff}
\begin{align}
 &h(d\phi(X),d\phi(Y)) = \lambda^{2} g(X,Y) ,\quad \forall X,Y \in \h
\label{eqa} ; \\
& - \Delta^{\phi}\tau(\phi) - \tr{{\mathrm R}^{N}(d\phi,
\tau(\phi))d\phi} = 0 \label{eqb};\\
&  \lambda^{2} \tau(\phi) + d\phi \g{\lambda^{2}} = 0 \label{eqc}; \\
& |\tau(\phi)|^{4} - 2\Delta \lambda^{2} |\tau(\phi)|^{2} + 4 \Delta \lambda^{2}
\di{\langle d\phi , \tau(\phi)\rangle} + n (\Delta \lambda^{2})^{2} \label{eqd}\\
&+ 2 \langle d\phi
, \tau(\phi) \rangle (\nabla |\tau(\phi)|^{2}) + |S|^{2} = 0  .\notag
\end{align}
Equip the manifold $(M,g)$ with harmonic coordinates $(x^{i})_{1 \leqslant i\leqslant
m}$ centered at the point $p$ and $(N,h)$ with harmonic coordinates
$(y^{\alpha})_{1\leqslant\alpha\leqslant n}$ around the point $\phi(p)$.
\newline Let $f : U \subset N \to {\rn}$ be a local function on $N$, then
\begin{align}
& \Delta^{2} (f\co\phi)  = \frac{\partial^{4} f}{\partial y^{\alpha}\partial
y^{\beta}\partial y^{\gamma}\partial y^{\delta}} \left[ g^{ij}g^{kl} \frac{\partial
\phi^{\alpha}}{\partial x^{i}} \frac{\partial \phi^{\beta}}{\partial x^{j}}
\frac{\partial \phi^{\gamma}}{\partial x^{k}} \frac{\partial
\phi^{\delta}}{\partial x^{l}} \right] + \label{eqe}\\
& \frac{\partial^{3} f}{\partial y^{\alpha}\partial y^{\beta}\partial y^{\gamma}}
\left[ \left( \frac{\partial^{2} \phi^{\gamma}}{\partial x^{i}\partial x^{j}}
\frac{\partial \phi^{\beta}}{\partial x^{k}}\frac{\partial \phi^{\alpha}}{\partial
x^{l}} + \frac{\partial^{2} \phi^{\beta}}{\partial x^{i}\partial x^{k}}
\frac{\partial \phi^{\gamma}}{\partial x^{j}}\frac{\partial \phi^{\alpha}}{\partial
x^{l}} + \frac{\partial^{2} \phi^{\alpha}}{\partial x^{i}\partial x^{l}}
\frac{\partial \phi^{\gamma}}{\partial x^{j}}\frac{\partial
\phi^{\beta}}{\partial x^{k}} + \right. \right. \notag\\
& \left. \frac{\partial^{2} \phi^{\beta}}{\partial x^{k}\partial x^{j}}
\frac{\partial \phi^{\gamma}}{\partial x^{i}}\frac{\partial \phi^{\alpha}}{\partial
x^{l}} + \frac{\partial^{2} \phi^{\alpha}}{\partial x^{l}\partial x^{j}}
\frac{\partial \phi^{\gamma}}{\partial x^{i}}\frac{\partial \phi^{\beta}}{\partial
x^{k}} + \frac{\partial^{2} \phi^{\alpha}}{\partial x^{l}\partial x^{k}}
\frac{\partial \phi^{\gamma}}{\partial x^{i}}\frac{\partial
\phi^{\beta}}{\partial x^{j}} \right) g^{ij}g^{kl} +  \notag\\
& \left. 2 g^{ij} \frac{\partial g^{kl}}{\partial x^{i}} \left( \frac{\partial
\phi^{\gamma}}{\partial x^{j}} \frac{\partial \phi^{\beta}}{\partial x^{k}}
\frac{\partial \phi^{\alpha}}{\partial x^{l}} \right) \right]  + \notag\\
& \frac{\partial^{2} f}{\partial y^{\alpha} \partial y^{\beta}} \left[ \left(
\frac{\partial^{3} \phi^{\beta}}{\partial x^{i}\partial x^{j}\partial x^{k}}
\frac{\partial \phi^{\alpha}}{\partial x^{l}} + \frac{\partial^{2}
\phi^{\beta}}{\partial x^{j}\partial x^{k}} \frac{\partial^{2}
\phi^{\alpha}}{\partial x^{l}\partial x^{i}} + \frac{\partial^{2}
\phi^{\beta}}{\partial x^{i}\partial x^{k}} \frac{\partial^{2}
\phi^{\alpha}}{\partial x^{l}\partial x^{j}} +
\right. \right. \notag\\
& \left. \frac{\partial^{3} \phi^{\alpha}}{\partial x^{i}\partial x^{j}\partial
x^{l}} \frac{\partial \phi^{\beta}}{\partial x^{k}} + \frac{\partial^{2}
\phi^{\beta}}{\partial x^{i}\partial x^{j}} \frac{\partial^{2}
\phi^{\alpha}}{\partial x^{l}\partial x^{k}} + \frac{\partial^{3}
\phi^{\alpha}}{\partial x^{i}\partial x^{k}\partial x^{l}} \frac{\partial
\phi^{\beta}}{\partial x^{j}} + \frac{\partial^{3} \phi^{\alpha}}{\partial
x^{j}\partial x^{k}\partial x^{l}}
\frac{\partial \phi^{\beta}}{\partial x^{i}} \right) g^{ij}g^{kl} + \notag\\
&\left. \left( \frac{\partial^{2} \phi^{\alpha}}{\partial x^{k}\partial x^{l}}
\frac{\partial \phi^{\beta}}{\partial x^{j}} + \frac{\partial^{2}
\phi^{\alpha}}{\partial x^{j}\partial x^{l}} \frac{\partial \phi^{\beta}}{\partial
x^{k}} + \frac{\partial^{2} \phi^{\beta}}{\partial x^{j}\partial x^{k}}
\frac{\partial \phi^{\alpha}}{\partial x^{l}} \right) 2 g^{ij} \frac{\partial
g^{kl}}{\partial x^{i}} + g^{ij} \frac{\partial^{2} g^{kl}}{\partial x^{i} \partial
x^{j}} \frac{\partial \phi^{\alpha}}{\partial x^{l}}
\frac{\partial \phi^{\beta}}{\partial x^{k}} \right] + \notag\\
&\frac{\partial f}{\partial y^{\alpha}} \left[ \frac{\partial^{2}
\phi^{\alpha}}{\partial x^{k}\partial x^{l}} g^{ij} \frac{\partial^{2}
g^{kl}}{\partial x^{i} \partial x^{j}} + \frac{\partial^{3} \phi^{\alpha}}{\partial
x^{j}\partial x^{k}\partial x^{l}} 2 g^{ij} \frac{\partial g^{kl}}{\partial x^{i}} +
\frac{\partial^{4} \phi^{\alpha}}{\partial x^{i}\partial x^{j}\partial x^{k}\partial
x^{l}}g^{ij}g^{kl} \right] . \notag
\end{align}

Plugging carefully chosen local biharmonic test functions, as given
by~\cite[Proposition~2.4]{A-G}, into Equation~\eqref{eqe} shows that $\phi$ is a
biharmonic morphism if and only if
\begin{align}
& g^{ij} \frac{\partial \phi^{\alpha}}{\partial x^{i}}\frac{\partial
\phi^{\beta}}{\partial x^{j}} = \lambda^{2} h^{\alpha\beta} \label{s1};\\
&d\phi^{\gamma} \g{\left(\lambda^{2}\right)} + \lambda^{2}
\tau^{\gamma}(\phi)  = 0 \label{s2};\\
&h^{\alpha\beta}\Delta(\lambda^{2}) + 
\di{ (\tau^{\beta}(\phi)d\phi^{\alpha})} +
\di{ (\tau^{\alpha}(\phi)d\phi^{\beta}}) -
\tau^{\alpha}(\phi)\tau^{\beta}(\phi) =\lambda^{2} dh^{\alpha\beta}(\tau(\phi)) \label{s3};\\
&g^{ij} \frac{\partial^{2} }{\partial x^{i} \partial x^{j}} \tau^{\alpha}(\phi) = 0,
\label{s4}
\end{align}
for all $\alpha, \beta = 1,\dots, n$.\\
Clearly Equations~\eqref{s1} and \eqref{s2} mean that $\phi$ is horizontally weakly conformal and $4$-harmonic.\\
On the other hand, Equation~\eqref{s3} is equivalent to the vanishing of
$$A^{\alpha\beta} = h^{\alpha\beta}\Delta(\lambda^{2}) +
g^{ij}\frac{\partial \phi^{\alpha}}{\partial x^{i}} \frac{\partial
\tau^{\beta}(\phi)}{\partial x^{j}} + g^{ij}\frac{\partial \phi^{\beta}}{\partial
x^{i}} \frac{\partial \tau^{\alpha}(\phi)}{\partial x^{j}}
+\tau^{\alpha}(\phi)\tau^{\beta}(\phi) - \lambda^{2} \frac{\partial
h^{\alpha\beta}}{\partial y^{p}} \tau^{p}(\phi),$$ whose norm is easily shown to be
\begin{align*}
|A|^{2} &=  n (\Delta \lambda^{2})^{2} + |\tau(\phi)|^{4} -
2(\Delta\lambda^{2})|\tau(\phi)|^{2} + 4 (\Delta\lambda^{2})
\di{\langle d\phi , \tau(\phi)\rangle} \label{eq10}\\
& + 4 g^{ij} \frac{\partial \phi^{\alpha}}{\partial x^{i}} \frac{\partial
\tau^{\beta}(\phi)}{\partial x^{j}} h_{\alpha\delta}h_{\beta\mu}
\tau^{\delta}(\phi)\tau^{\mu}(\phi)
  + 2 \lambda^{2} \tau^{\alpha}(\phi)\tau^{\beta}(\phi)
\frac{\partial h_{\alpha\beta}}{\partial y^{p}} \tau^{p}(\phi) \notag\\
& + 2 g^{ij} \frac{\partial \phi^{\alpha}}{\partial x^{i}} \frac{\partial
\tau^{\beta}(\phi)}{\partial x^{j}} h_{\alpha\delta}h_{\beta\mu} g^{kl}
\frac{\partial \phi^{\delta}}{\partial x^{k}} \frac{\partial
\tau^{\mu}(\phi)}{\partial x^{l}} +  2 g^{ij} \frac{\partial \phi^{\alpha}}{\partial
x^{i}} \frac{\partial \tau^{\beta}(\phi)}{\partial x^{j}}
h_{\alpha\delta}h_{\beta\mu} g^{kl} \frac{\partial \phi^{\mu}}{\partial x^{k}}
\frac{\partial
\tau^{\delta}(\phi)}{\partial x^{l}} \notag\\
& + 4 \lambda^{2} g^{ij} \frac{\partial \phi^{\alpha}}{\partial x^{i}} \frac{\partial
\tau^{\beta}(\phi)}{\partial x^{j}} \frac{\partial h_{\alpha\beta}}{\partial y^{p}}
\tau^{p}(\phi) + \lambda^{4} \frac{\partial h^{\alpha\beta}}{\partial y^{p}}
h_{\alpha\delta}h_{\beta\mu} \frac{\partial h^{\delta\mu}}{\partial y^{q}}
\tau^{p}(\phi) \tau^{q}(\phi) , \notag
\end{align*}
since $h^{\alpha\beta}h_{\alpha\delta}h_{\beta\mu} = h_{\delta\mu}$ and
$h_{\alpha\delta}h_{\beta\mu}\frac{\partial h^{\alpha\beta}}{\partial y^{p}} = -
\frac{\partial
h_{\alpha\beta}}{\partial y^{p}}$.\\
But
\begin{align*}
\langle d\phi, \tau(\phi) \rangle (\nabla |\tau(\phi)|^{2}) = 2 g^{ij} \frac{\partial
\phi^{\alpha}}{\partial x^{i}} \frac{\partial \tau^{\beta}(\phi)}{\partial x^{j}}
h_{\alpha\delta}h_{\beta\mu} \tau^{\delta}(\phi)\tau^{\mu}(\phi)
  +  \lambda^{2} \tau^{\alpha}(\phi)\tau^{\beta}(\phi)
\frac{\partial h_{\alpha\beta}}{\partial y^{p}} \tau^{p}(\phi)
\end{align*}
and the components of $S \in \odot^{2} \phi^{-1}TN$, the symmetrization of $g(d\phi ,
\nabla^{\phi}\tau(\phi) )$, are
$$ S^{\alpha\beta} =  g^{ij} \frac{\partial \phi^{\alpha}}{\partial x^{i}}
\frac{\partial \tau^{\beta}(\phi)}{\partial x^{j}} + \lambda^{2} h^{\alpha\delta}
\Gamma^{\beta}_{\gamma\delta} \tau^{\gamma}(\phi)
  + g^{ij} \frac{\partial \phi^{\beta}}{\partial x^{i}}
\frac{\partial \tau^{\alpha}(\phi)}{\partial x^{j}} + \lambda^{2} h^{\beta\delta}
\Gamma^{\alpha}_{\gamma\delta} \tau^{\gamma}(\phi) .
$$
So its norm is
\begin{align*}
|S|^{2} &= 2 g^{ij} \frac{\partial \phi^{\alpha}}{\partial x^{i}} \frac{\partial
\tau^{\beta}(\phi)}{\partial x^{j}} h_{\alpha\delta}h_{\beta\mu} g^{kl}
\frac{\partial \phi^{\delta}}{\partial x^{k}} \frac{\partial
\tau^{\mu}(\phi)}{\partial x^{l}} +  2 g^{ij} \frac{\partial \phi^{\alpha}}{\partial
x^{i}} \frac{\partial \tau^{\beta}(\phi)}{\partial x^{j}}
h_{\alpha\delta}h_{\beta\mu} g^{kl} \frac{\partial \phi^{\mu}}{\partial x^{k}}
\frac{\partial
\tau^{\delta}(\phi)}{\partial x^{l}} \\
& + 4 \lambda^{2} g^{ij} \frac{\partial \phi^{\alpha}}{\partial x^{i}} \frac{\partial
\tau^{\beta}(\phi)}{\partial x^{j}} \frac{\partial h_{\alpha\beta}}{\partial y^{p}}
\tau^{p}(\phi) + \lambda^{4} \frac{\partial h^{\alpha\beta}}{\partial y^{p}}
h_{\alpha\delta}h_{\beta\mu} \frac{\partial h^{\delta\mu}}{\partial y^{q}}
\tau^{p}(\phi) \tau^{q}(\phi).
\end{align*}
Hence \eqref{s3} is equivalent to Equation~\eqref{longeq}.\\
To obtain that Equation~\eqref{s4} is biharmonicity, observe that
\begin{align*}
&- (\Delta^{\phi}  \tau(\phi) )^{\alpha}= g^{ij} \left[ \frac{\partial^{2}
\tau^{\alpha}(\phi)}{\partial x^{i}\partial x^{j}}  + \frac{\partial
\tau^{\gamma}(\phi)}{\partial x^{j}} \frac{\partial \phi^{\beta}}{\partial x^{i}}
\Gamma^{\alpha}_{\beta\gamma} + \frac{\partial \tau^{\gamma}(\phi)}{\partial x^{i}}
\frac{\partial \phi^{\beta}}{\partial x^{j}} \Gamma^{\alpha}_{\beta\gamma} +
\tau^{\gamma}(\phi) \frac{\partial^{2} \phi^{\beta}}{\partial
x^{i}\partial x^{j}} \Gamma^{\alpha}_{\beta\gamma} \right.\\
  &+ \left.
\tau^{\gamma}(\phi) \frac{\partial \phi^{\beta}}{\partial x^{j}} \frac{\partial
\Gamma^{\alpha}_{\beta\gamma}}{\partial y^{\mu}} \frac{\partial \phi^{\mu}}{\partial
x^{i}} + \tau^{\mu}(\phi) \frac{\partial \phi^{\beta}}{\partial x^{j}}
\Gamma^{\gamma}_{\mu\beta} \frac{\partial \phi^{\delta}}{\partial x^{i}}
\Gamma^{\alpha}_{\gamma\delta} - \Gamma^{k}_{ij} \frac{\partial
\tau^{\alpha}(\phi)}{\partial x^{k}}
  - \Gamma^{k}_{ij} \tau^{\gamma}(\phi)
\frac{\partial \phi^{\beta}}{\partial x^{k}}
\Gamma^{\alpha}_{\beta\gamma} \right] \\
& = 2 g^{ij} \frac{\partial \tau^{\gamma}(\phi)}{\partial x^{i}} \frac{\partial
\phi^{\beta}}{\partial x^{j}} \Gamma^{\alpha}_{\beta\gamma} + \tau^{\gamma}(\phi)
\tau^{\beta}(\phi) \Gamma^{\alpha}_{\beta\gamma} + \lambda^{2} \tau^{\gamma}(\phi)
h^{\beta\mu} \frac{\partial \Gamma^{\alpha}_{\beta\gamma}}{\partial y^{\mu}} +
\lambda^{2} h^{\beta\delta}\tau^{\mu}(\phi) \Gamma^{\gamma}_{\mu\beta}
\Gamma^{\alpha}_{\gamma\delta} ,
\end{align*}
since $\phi$ is {\hwc}, $g^{ij} \frac{\partial^{2} }{\partial x^{i} \partial x^{j}}
\tau^{\alpha_{0}}(\phi)=0$ and we are working with harmonic coordinates.\\
On the other hand, $\tr{{\mathrm R}^{N}(d\phi, \tau(\phi))d\phi} = \lambda^{2}
{\mathrm {Ric}}^{N} \tau(\phi)$ and
\begin{align*}
- \lambda^{2} {\mathrm {Ric}}^{N} \tau(\phi) &=  \lambda^{2} \left[ - \frac{\partial
h^{\alpha\beta}}{\partial y^{\gamma}} \Gamma^{\delta}_{\alpha\beta} - h^{\alpha\beta}
\frac{\partial \Gamma^{\delta}_{\alpha\gamma}}{\partial y^{\beta}} - \sum_{\mu
=1}^{n} h^{\alpha\beta} \Gamma^{\mu}_{\alpha\gamma} \Gamma^{\delta}_{\mu\beta}
\right] \tau^{\gamma}(\phi) .
\end{align*}
Thus
\begin{align*}
\tau^{2}(\phi)  &= \left[ 2 g^{ij} \frac{\partial \tau^{\gamma}(\phi)}{\partial
x^{i}} \frac{\partial \phi^{\beta}}{\partial x^{j}} \Gamma^{\alpha}_{\beta\gamma} +
\tau^{\gamma}(\phi) \tau^{\beta}(\phi) \Gamma^{\alpha}_{\beta\gamma} - \lambda^{2}
\frac{\partial h^{\delta\beta}}{\partial y^{\gamma}} \Gamma^{\alpha}_{\delta\beta}
\tau^{\gamma}(\phi)\right]  \frac{\partial }{\partial y^{\alpha}}= 0.
\end{align*}
Of course, Conditions~\eqref{s1}--\eqref{s4} are also sufficient.
\end{proof}

It is well-known that harmonic morphisms into $2$-dimensional manifolds have an
interesting link to the geometry of the fibres, which can be stated as: a horizontally
conformal submersion $(M^m, g)\longrightarrow (N^2,h)$ is harmonic (hence a harmonic
morphism) if and only if it has minimal fibres. The corresponding result for
biharmonic morphisms is partially true.
\begin{cor}
A biharmonic morphism $\phi : (M^m,g)\longrightarrow (N^4,h)$ ($m\geqslant 4$), with
a $4$-dimensional target, always has minimal fibres. However there exist horizontally
conformal submersions which are $4$-harmonic, have minimal fibres but are not
biharmonic morphisms.
\end{cor}
\begin{proof}
By Theorem~\ref{thm1}, a biharmonic morphism is a $4$-harmonic morphism. The first
statement then follows from a result in~\cite{B-G} that a horizontally weakly
conformal map $(M^m, g)\longrightarrow (N^n,h)$
is a $p$-harmonic map (hence a $p$-harmonic morphism) with $p=n=\dim{N}$ if and only if it has minimal fibres.\\
It is well-known that the radial projection $\varphi: \mathbb{R}^{5}\setminus
\{0\}\longrightarrow {\sn}^{4}, \; \varphi(x)=x/|x|$ is a horizontally homothetic
submersion with totally geodesic fibres and hence a harmonic morphism
(see~\cite{B-W}). It is also a $p$-harmonic morphism for any $p$ and, in particular,
a $4$-harmonic morphism. However, the radial projection $\varphi:
\mathbb{R}^{m}\setminus \{0\}\longrightarrow {\sn}^{m-1}, \; \varphi(x)=x/|x|$ is a
biharmonic morphism if and only if $m=4$ (\cite{Ou}).
\end{proof}

\begin{rk}
1) A straightforward computation shows that if $\phi$ is a biharmonic morphism then
there exists a continuous function $\lambda$ on $M$ such that
\begin{equation*}
\tau^{2} (\psi\circ\phi) = \lambda^{4} \tau^{2}(\psi) \circ \phi ,
\end{equation*}
for any map $\psi$ (including functions).\\
2) Taking the trace of Equation~\eqref{s3} yields
\begin{equation}\label{trace}
n \Delta(\lambda^{2}) + 2 \di{\langle d\phi , \tau(\phi)\rangle} = |\tau(\phi)|^{2} ,
\end{equation}
hence by Stokes' Theorem, biharmonic morphisms on a compact manifold without boundary
 are exactly homothetic submersions with minimal fibres. See~\cite{L-O} for similar results.
\end{rk}

Apart from harmonic morphisms with harmonic dilation, examples of biharmonic
morphisms are hard to unearth. Conformal maps in dimension four have the double
advantage of satisfying automatically two of the four conditions, namely horizontal
weak conformality and $4$-harmonicity. Once biharmonicity is secured, remains only
Equation~\eqref{longeq} to fulfil, though no geometric insight is yet at our
disposal.

\begin{thm}
The inversion in the unit sphere
\begin{align*}
\phi : \rn^{n}\setminus\{0\} &\to \rn^n \\
x &\mapsto \frac{x}{|x|^2}
\end{align*}
is a biharmonic morphism if and only if $n=4$.
\end{thm}

\begin{proof}
By \cite{B-K}, the inversion in the unit sphere is a biharmonic map if and only if
$n=4$. Moreover the map
\begin{align*}
\phi : \rn^{4}\setminus\{0\} &\to \rn^4 \\
x &\mapsto \frac{x}{|x|^2}
\end{align*}
is clearly a conformal map of dilation $\lambda^2 = \frac{1}{|x|^4}$ between spaces
of dimension four, hence a $4$-harmonic map (\cite{Ou}) and $\tau(\phi) = -4
\frac{x}{|x|^4}$. Using the standard coordinates $\{ x^{\alpha} \}_{\alpha =
1,2,3,4}$ on $\rn^4$, simple computations show that
\begin{align*}
|\tau(\phi)|^2 &= \frac{16}{|x|^6} ; \quad \Delta \lambda^2 = \frac{8}{|x|^6} ; \quad
\langle d\phi,\tau(\phi)\rangle = 4\frac{x}{|x|^6} ; \quad \di \langle
d\phi,\tau(\phi)\rangle = -\frac{8}{|x|^6}.
\end{align*}
As to the symmetric tensor $S$, we have
\begin{align*}
S^{\alpha\alpha} &= -\tfrac{8}{|x|^8} \left( |x|^2 + 2 (x^{\alpha})^2 \right) ; \quad
S^{\alpha\beta} = -\tfrac{16}{|x|^8}  x^{\alpha}x^{\beta} ,
\end{align*}
for all $\alpha \neq \beta = 1,\dots,4$, so $|S|^2 = \tfrac{3 (16)^2}{|x|^{12}}$.
Therefore $\phi$ satisfies Equation~\eqref{longeq} and is a biharmonic morphism.
\end{proof}

\section{Biharmonicity and conformality}

Since biharmonic morphisms appear so rigid, we drop two of their characteristic
conditions and only keep horizontal weak conformality and biharmonicity. In the light of
the Fuglede-Ishihara theorem, such maps are the exact counterparts of harmonic
morphisms. Moreover, the expression of the tension field of a horizontally weakly
conformal map enlightens the relationship between minimality of the fibres and
harmonicity of the map. Though the formula for the bitension field is far more
intricate, we can apply it to some special cases to obtain new examples of
biharmonic maps.

We will denote by $\mu = \tfrac{1}{m-n} \sum_{s=1}^{m-n} (\nabla_{e_{s}} e_{s})^\h$
and $\nu= \tfrac{1}{n}\sum_{i=1}^{n} (\nabla_{e_{i}} e_{i})^\v$, for an orthonormal
frame $\{e_{i},e_{s}\}_{i=1,\dots,n, s=1,\dots,m-n}$, with $e_{i}$ horizontal and
$e_{s}$ vertical, the mean curvatures of the vertical and horizontal distributions,
and $A$ and $B$ the second fundamental forms of the horizontal and vertical
distributions:
$$ A_{E}F = (\nabla_{E^\h} F^\h)^\v \quad B_{E}F = (\nabla_{E^\v}F^\v)^\h, \quad E,F \in \Gamma(TM).$$

\begin{thm}\label{prop1}
Let $\phi : (M^m ,g) \to (N^n ,h)$ ($m\geqslant n \geqslant 2$) be a horizontally
weakly conformal map of dilation $\lambda$ between Riemannian manifolds, then
$\phi$ is biharmonic if and only if, at every regular point
\begin{align*}
& d\phi(-\Delta X) + \Delta^\h (\ln\lambda) d\phi(X) + (n-2) \nu(\ln\lambda)d\phi(X)
  + 2 d\phi(\nabla_{\g^\h \ln\lambda} X) \\
  &+ 2d\phi(\g X(\ln\lambda))-d\phi(\nabla_{X} (\g^\h\ln\lambda)) -2 (\di^\h X) d\phi(\g \ln\lambda)\\
& + X(\ln \lambda)d\phi(X) 
 - d\phi (\nabla_{X} \g \ln\lambda) + \lambda^2 \ricci^{N}(d\phi(X))
 + (n-1) d\phi(\nabla_{X} \nu)  \\
 &- (m-n)[ \mu(\ln\lambda) d\phi(X) - \langle X, \mu \rangle d\phi(\g\ln\lambda)]
-(m-n)d\phi (\nabla_{X} \mu)+ n d\phi(A^{*}_{X}\nu) \\
&+ \tr d\phi( (\nabla A)^{*}_{X} - (\nabla A)X + 3A^{*}_{\nabla X} + (\nabla_{X}
B^{*})^{*} - B_{\nabla X} + 2 B_{\nabla_{X}} -2 A^{*}_{\nabla_{X}} ) = 0,
\end{align*}
where $X= (2-n)\g^{\h}\ln \lambda - (m-n)\mu$, $\g^{\h}$ being the horizontal
gradient.
\end{thm}
\begin{proof}
At a regular point, the tension field of a horizontally weakly conformal map $\phi :
(M^m ,g) \to (N^n ,h)$, of dilation $\lambda$, is $\tau(\phi) = d\phi(X)$, where $X=
(2-n)\g^{\h}\ln \lambda - (m-n)\mu$. Recall that the mean curvature of the horizontal
distribution is $\nu = \g^\v \ln\lambda$. Then, for an adapted orthonormal frame
$\{e_{i}, e_{s}\}_{i=1,\dots,n, s=1,\dots,m-n}$, $e_{i} \in \h$ and $e_{s} \in \v$ on
$(M^m ,g)$
\begin{align*}
\nabla^{\phi}_{e_{i}} d\phi(X)
&= e_{i}(\ln \lambda) d\phi(X) + X(\ln \lambda) d\phi(e_{i}) - \langle X,e_{i}\rangle
d\phi(\g \ln \lambda) + d\phi(\nabla_{e_{i}}X) ,
\end{align*}
so
\begin{align*}
&\nabla^{\phi}_{e_{i}} \nabla^{\phi}_{e_{i}} d\phi(X) 
= e_{i}(e_{i}(\ln \lambda)) d\phi(X) + (e_{i}(\ln \lambda))^2 d\phi(X) + e_{i}(\ln \lambda)X(\ln \lambda)d\phi(e_{i}) \\
&-e_{i}(\ln \lambda)\langle X,e_{i} \rangle d\phi(\g\ln \lambda)
+ e_{i}(\ln \lambda)d\phi(\nabla_{e_{i}}X) + e_{i}(X(\ln\lambda))d\phi(e_{i}) \\
&+ X(\ln\lambda) \nabla^{\phi}_{e_{i}} d\phi(e_{i}) - e_{i}\langle X,e_{i} \rangle
d\phi(\g\ln \lambda)
- \langle X,e_{i} \rangle e_{i}(\ln \lambda) d\phi(\g\ln \lambda) \\
&- \langle X,e_{i} \rangle (\g^\h\ln \lambda) (\ln \lambda) d\phi(e_{i})
+ \langle X,e_{i} \rangle \langle \g\ln\lambda ,e_{i} \rangle d\phi(\g\ln \lambda) \\
&- \langle X,e_{i} \rangle d\phi(\nabla_{e_{i}} \g^\h\ln\lambda) + e_{i}(\ln \lambda)d\phi(\nabla_{e_{i}}X) \\
&+ (\nabla_{e_{i}}X)^\h(\ln\lambda)d\phi(e_{i}) - \langle e_{i},\nabla_{e_{i}}X
\rangle d\phi(\g\ln\lambda) + d\phi(\nabla_{e_{i}}(\nabla_{e_{i}}X)^\h) .
\end{align*}
Summing on the index $i$
\begin{align*}
&\big( \nabla^{\phi}_{e_{i}} \nabla^{\phi}_{e_{i}} -
\nabla^{\phi}_{\nabla_{e_{i}}e_{i}}\big) d\phi(X) =
e_{i}(e_{i}(\ln \lambda)) d\phi(X) + 2 d\phi(\nabla_{\g^\h \ln\lambda} X) \\
&+ d\phi(\g (X(\ln\lambda))) 
+ X(\ln\lambda)(\nabla d\phi)(e_{i},e_{i}) -d\phi(\nabla_{X} (\g^\h\ln\lambda)) \\
&+ (\nabla_{e_{i}}X)^\h (\ln\lambda)d\phi(e_{i}) 
- 2\langle e_{i},\nabla_{e_{i}}X \rangle d\phi(\g\ln\lambda)
 + d\phi(\nabla_{e_{i}}(\nabla_{e_{i}}X)^\h)\\
&- (\nabla_{e_{i}}e_{i})^\h (\ln \lambda) d\phi(X) +
d\phi(\nabla_{X}(\nabla_{e_{i}}e_{i})^\v ) - d\phi(\nabla_{\nabla_{e_{i}}e_{i}}X).
\end{align*}
On the other hand, for the vertical bundle
\begin{align*}
- \nabla^{\phi}_{\nabla_{e_{s}}e_{s}} d\phi(X) &= - (\nabla_{e_{s}}e_{s} )^\h (\ln
\lambda) d\phi(X)- X(\ln \lambda) d\phi(\nabla_{e_{s}}e_{s})\\
&+ \langle X, \nabla_{e_{s}}e_{s} \rangle d\phi(\g\ln\lambda)  
- d\phi ( \nabla_{(\nabla_{e_{s}}e_{s})^\h} X ) - d\phi
([(\nabla_{e_{s}}e_{s})^\v,X]),
\end{align*}
so
\begin{align*}
&(\nabla^{\phi}_{e_{s}} \nabla^{\phi}_{e_{s}} - \nabla^{\phi}_{\nabla_{e_{s}}e_{s}})
d\phi(X) = d\phi(\nabla_{e_{s}}\nabla_{e_{s}} X) - d\phi(\nabla_{e_{s}}\nabla_{X}
e_{s})  - d\phi(\nabla_{\nabla_{e_{s}}X}e_{s})\\
&+ d\phi(\nabla_{\nabla_{X}e_{s}}e_{s}) 
- (\nabla_{e_{s}}e_{s} )^\h (\ln \lambda) d\phi(X) - X(\ln \lambda)
d\phi(\nabla_{e_{s}}e_{s})\\
& + \langle X, \nabla_{e_{s}}e_{s} \rangle
d\phi(\g\ln\lambda)- d\phi (\nabla_{\nabla_{e_{s}}e_{s}}X) 
+ d\phi (\nabla_{X}(\nabla_{e_{s}}e_{s})^\v).
\end{align*}
Therefore
\begin{align*}
&-\Delta (\tau(\phi))
= d\phi(-\Delta X) - d\phi(\nabla_{e_{i}}(\nabla_{e_{i}}X)^\v) + \Delta^\h (\ln\lambda) d\phi(X) + n\nu(\ln\lambda) d\phi(X) \\
& + 2 d\phi(\nabla_{\g^\h \ln\lambda} X) + d\phi(\g X(\ln\lambda)) -d\phi(\nabla_{X} (\g^\h\ln\lambda)) -2 (\di^\h X) d\phi(\g \ln\lambda) \\
& + X(\ln \lambda)(\tau(\phi) + (m-n) d\phi(\mu) )
+ (\nabla_{e_{i}}X)^\h (\ln\lambda)d\phi(e_{i})  + n d\phi(\nabla_{X} \nu) \\
& - (m-n)[ \mu(\ln\lambda) d\phi(X) + X(\ln\lambda) d\phi(\mu) - \langle X, \mu \rangle d\phi(\g\ln\lambda)] \\
& - d\phi(\nabla_{e_{s}}\nabla_{X} e_{s}) - d\phi(\nabla_{\nabla_{e_{s}}X}e_{s}) +
d\phi(\nabla_{\nabla_{X}e_{s}}e_{s})+ d\phi (\nabla_{X}(\nabla_{e_{s}}e_{s})^\v ).
\end{align*}
But
\begin{align*}
&(\nabla_{e_{i}}X)^\h (\ln\lambda)d\phi(e_{i}) 
= d\phi( \g  X (\ln\lambda ) ) - d\phi(\nabla_{X} \g \ln\lambda)\\
& - \langle (\nabla_{e_{i}}X)^\v ,  \g^\v \ln\lambda \rangle d\phi(e_{i}).
\end{align*}
However (cf.~\cite[Proposition 2.5.17]{B-W})
\begin{align*}
- \langle (\nabla_{e_{i}}X)^\v ,  \g^\v \ln\lambda \rangle d\phi(e_{i}) &=
 - \langle \tfrac{1}{2} (\nabla_{e_{i}}X)^\v ,  \g^\v \ln\lambda \rangle d\phi(e_{i}) \\
&+ \langle \tfrac{1}{2} (\nabla_{X}e_{i})^\v ,  \g^\v \ln\lambda \rangle d\phi(e_{i})
- |\g^\v \ln\lambda|^2 d\phi(X)
\end{align*}
since $ ({\mathcal L}_{V} g)(X,Y) = - g( \nabla_{X}Y + \nabla_{Y}X , V) = -d(\ln\lambda^2) (V) g(X,Y)$.
Hence
\begin{align*}
&- \langle (\nabla_{e_{i}}X)^\v ,  \g^\v \ln\lambda \rangle d\phi(e_{i}) = - d\phi(\nabla_{X}\g^\v \ln\lambda) - 2|\g^\v \ln\lambda|^2 d\phi(X).
\end{align*}
Therefore
\begin{align*}
(\nabla_{e_{i}}X)^\h (\ln\lambda)d\phi(e_{i}) &=  d\phi( \g  X (\ln\lambda ) ) - d\phi(\nabla_{X} \g \ln\lambda)- d\phi(\nabla_{X}\nu) \\
&- 2|\nu|^2 d\phi(X)
\end{align*}
and
\begin{align*}
&-\Delta (\tau(\phi)) =  d\phi(-\Delta X) - d\phi(\nabla_{e_{i}}(\nabla_{e_{i}}X)^\v)
+ \Delta^\h (\ln\lambda) d\phi(X) \\
&+ (n-2) \nu(\ln\lambda)d\phi(X) 
 + 2 d\phi(\nabla_{\g^\h \ln\lambda} X) + 2d\phi(\g X(\ln\lambda)) \\
&-d\phi(\nabla_{X} (\g^\h\ln\lambda)) -2 (\di^\h X) d\phi(\g \ln\lambda)
+ X(\ln \lambda)\tau(\phi) - d\phi (\nabla_{X} \g \ln\lambda) \\
& + (n-1) d\phi(\nabla_{X} \nu)  - (m-n)[ \mu(\ln\lambda) d\phi(X) - \langle X, \mu \rangle d\phi(\g\ln\lambda)] \\
& - d\phi(\nabla_{e_{s}}\nabla_{X} e_{s}) - d\phi(\nabla_{\nabla_{e_{s}}X}e_{s}) +
d\phi(\nabla_{\nabla_{X}e_{s}}e_{s})+ d\phi (\nabla_{X}(\nabla_{e_{s}}e_{s})^\v ).
\end{align*}

\begin{lem}
Let $A$ and $B$ be the second fundamental forms of the horizontal and vertical
distributions. Then
\begin{align*}
&d\phi(\tr_{h} \nabla (\nabla X)^\v) = d\phi(\tr (\nabla A)X) ;\\
&d\phi(-\nabla_{e_{s}}\nabla_{X} e_{s} - \nabla_{\nabla_{e_{s}}X}e_{s}
+ \nabla_{\nabla_{X}e_{s}}e_{s} + \nabla_{X}(\nabla_{e_{s}}e_{s})^\v) = -(m-n)d\phi(\nabla_{X}\mu) \\
&+ \tr d\phi( (\nabla A)^{*}_{X} + 3 A^{*}_{\nabla X} + (\nabla_{X} B^{*})^{*} -
B_{\nabla X} -2A^{*}_{\nabla_{X}}  + 2 B_{\nabla_{X}}) + n d\phi(A^{*}_{X}\nu).
\end{align*}
\end{lem}

\begin{proof}
Consider an adapted orthonormal frame $\{e_{i}, e_{s}\}_{i=1,\dots,n,
s=1,\dots,m-n}$, $e_{i} \in \h$ and $e_{s} \in \v$. First observe that
$$ d\phi((\nabla_{e_{i}}A )_{e_{i}} X) =  d\phi(\nabla_{e_{i}} (\nabla_{e_{i}} X)^\v ),$$
and $(\nabla_{e_{s}} A)_{e_{s}} X$ is vertical.\\
For the second equality, we have
\begin{align*}
&\langle -\nabla_{e_{s}}\nabla_{X} e_{s} - \nabla_{\nabla_{e_{s}}X}e_{s}
+ \nabla_{\nabla_{X}e_{s}}e_{s} + \nabla_{X}\nabla_{e_{s}}e_{s}, e_{j}\rangle  \\
&= \langle (\nabla_{e_{s}}A)_{X} e_{j}, e_{s}\rangle  + \langle A_{\nabla_{e_{s}} X}
e_{j}, e_{s}\rangle + \langle (\nabla_{X}B^{*})_{e_{s}} e_{j}, e_{s}\rangle + \langle
B^{*}_{\nabla_{e_{s}}X} e_{j}, e_{s}\rangle \\
&+ 2 \langle \nabla_{[e_{s},X]} e_{j} ,e_{s} \rangle
\end{align*}
(confer~\cite[Th. 11.2.1 iii)]{B-W}.)\\
Moreover
\begin{align*}
\langle \nabla_{[e_{s},X]} e_{j} , e_{s} \rangle &= \langle e_{j} , A^{*}_{\nabla_{e_{s}}X} e_{s} -A^{*}_{\nabla_{X}e_{s}} e_{s} -
B_{\nabla_{e_{s}}X} e_{s} + B_{\nabla_{X}e_{s}} e_{s} \rangle .
\end{align*}
To conclude, observe that
\begin{align*}
 \langle (\nabla_{e_{i}}A)^{*}_{X} e_{i}, e_{j}\rangle &= - n\langle  e_{j} , A^{*}_{X}\nu \rangle ;\\
 A^{*}_{\nabla_{e_{i}} X} e_{i} &=  A^{*}_{\nabla_{X}e_{i}} e_{i} = 0 ;\\
\langle (\nabla_{X}B^{*})^{*}_{e_{i}} e_{i}, e_{j}\rangle &= 
 B_{\nabla_{e_{i}}X}e_{i} = B_{\nabla_{X} e_{i}}  e_{i}  = 0.
\end{align*}
\end{proof}
Note that for a horizontally weakly conformal map, the curvature term of the bitension field becomes
$$-\lambda^2 \ricci^{N}(d\phi(X)).$$
This completes the proof.
\end{proof}

Three special cases are interesting enough to be stated separately.

\begin{cor}\label{cor1}
i) A conformal map $\phi : (M^m,g) \to (N^n ,h)$ of conformal factor $\lambda$,
between manifolds of equal dimensions ($m=n > 2$), is biharmonic if and only if
\begin{align}
&- d\phi( \Delta (\g \ln\lambda)) -\Delta(\ln \lambda) d\phi(\g \ln \lambda) +  2 d\phi(\g |\g \ln \lambda|^2 )  \notag\\
& \quad + (2-n) |\g \ln \lambda|^2 d\phi(\g \ln \lambda)  + \lambda^2
\ricci^{N}(d\phi(\g\ln\lambda)) =0. \label{eqf}
\end{align}
ii) A horizontally conformal map $\phi : (M^m,g) \to (N^2 ,h)$ of dilation $\lambda$, into
a surface ($m > n=2$), is biharmonic if and only if
\begin{align*}
&d\phi(-\Delta \mu)+ \Delta^\h (\ln\lambda) d\phi(\mu) + 2 d\phi(\nabla_{\g^\h \ln\lambda} \mu) + 2d\phi(\g \mu(\ln\lambda)) \\
& -d\phi(\nabla_{\mu} (\g^\h\ln\lambda)) -2 (\di^\h \mu) d\phi(\g \ln\lambda)  \\
&   +(m-2) |\mu |^2 d\phi(\g \ln\lambda)  - d\phi (\nabla_{\mu} \g \ln\lambda) + \lambda^2\ricci^{N}(d\phi(\mu)) \\
&   + d\phi(\nabla_{\mu} \nu) + (2-m) d\phi(\nabla_{\mu}\mu)  -2(m-2)\mu(\ln\lambda) d\phi(\mu)\\
&+  \tr d\phi( (\nabla A)^{*}_{\mu} - (\nabla A)\mu + 3A^{*}_{\nabla \mu} +
(\nabla_{\mu} B^{*})^{*} + 2 B_{\nabla_{\mu}}- B_{\nabla \mu} -2
A^{*}_{\nabla_{\mu}}) \\
&+ 2 d\phi(A^{*}_{\mu}\nu) = 0.
\end{align*}
iii) A homothetic submersion $\phi : (M^m,g) \to (N^n ,h)$ ($m>n$), i.e. $\lambda$ is
constant, is biharmonic if and only if
\begin{align*}
&d\phi(\Delta \mu) - \lambda^2 \ricci^{N}(d\phi(\mu)) + (m-n) d\phi(\nabla_{\mu}\mu) \\
&- (\tr d\phi( (\nabla A)^{*}_{\mu} - (\nabla A)\mu + 3A^{*}_{\nabla \mu}
 + (\nabla_{\mu} B^{*})^{*} + 2 B_{\nabla_{\mu}}- B_{\nabla \mu} - 2 A^{*}_{\nabla_{\mu}} )) = 0.
\end{align*}
\end{cor}

\begin{rk} i) The case $m=n=2$ is trivial since any conformal map between surfaces is harmonic, hence biharmonic.\\
ii) When $\phi: (M^m ,g) \to (N^m ,h)$ is the identity and $h= e^{2\rho} g$,
Equation~\eqref{eqf} becomes Equation~(3.1) of~\cite{Bal}:
\begin{align*}
0 &= \tr_{g} \nabla^2 \g\rho + (-2\Delta \rho + (2-m)|\g\rho|^2) \g\rho \\
&+ \tfrac{6-m}{2} \g (|\g\rho|^2 ) + \ricci^{g}(\g\rho) .
\end{align*}
\end{rk}

The first case of Corollary~\ref{cor1} is the most prolific and the following
examples indicate that not only horizontally weakly conformal biharmonic maps are
easier to come by than biharmonic morphisms, but also that, even for this situation,
dimension four emerges as remarkable.

\begin{prop}
The inverse stereographic projection $\sigma_{\mathrm N}^{-1}$ from $(\rn^n ,ds^2)$
into $(\sn^{n}\setminus\{\mathrm N\} , ds^2)$, where $\mathrm N$ is the north pole and 
$ds^2$ the Euclidean metric, is a biharmonic map if and only if $n=4$.\\
Similarly, the identity map from  $(\bn^{n}, ds^2)$
into $(\bn^{n},\tfrac{4}{(1 - |x|^2)^2} ds^2)$, where $\bn^n$ is the open unit ball of $\rn^n$,
 is a biharmonic map if and only if $n=4$.\\
Furthermore, in either case, the biharmonic map is not a biharmonic morphism.
\end{prop}

\begin{proof}
For the first case, the map is given by $\sigma_{\mathrm N}^{-1}(x) = \tfrac{1}{1+ |x|^2}
(1-|x|^2 , 2x)$ and is isometric to the identity map from the Euclidean space $(\rn^n ,ds^2)$ into
($\rn^{n},\tfrac{4}{(1+ |x|^2)^2}ds^2)$. We can treat both examples at once, by considering the identity map 
$\phi$ on $\rn^n$ or $\bn^n$, from the Euclidean metric $ds^2$ into the conformal metric $\tfrac{4}{(1+ \epsilon |x|^2)^2}ds^2$
($\epsilon = 1$ for the sphere and $\epsilon = -1$ for the ball).\\
The map $\phi$ is  clearly conformal of dilation $\lambda = \tfrac{2}{1+ \epsilon |x|^2}$, so
\begin{align*}
\g \ln\lambda &= -\tfrac{2\epsilon}{1+ \epsilon |x|^2} x ; \quad |\g \ln\lambda|^2 = \tfrac{4|x|^2}{(1+ \epsilon |x|^2)^2}; \\
\g|\g \ln\lambda|^2 &= \tfrac{8(1- \epsilon |x|^2)}{(1+ \epsilon |x|^2)^3} x ; \quad
\Delta \ln\lambda = -\tfrac{2\epsilon}{(1+ \epsilon |x|^2)^2}(n + (n-2)\epsilon |x|^2) ;\\
(\tr\nabla^2) (\g\ln\lambda)  &= \tfrac{4}{(1+ \epsilon |x|^2)^3}(n+2 +
(n-2)\epsilon|x|^2)x,
\end{align*}
and Equation~\eqref{eqf} becomes
$$-8(n-4)\tfrac{(1- \epsilon|x|^2)}{(1+ \epsilon |x|^2)^3} x =0 ,$$
so both maps, into the sphere or the Poincar{\'e} model, are biharmonic if and only if $n=4$.\\
To know whether $\phi$ is also a biharmonic morphism we only
need to check Equation~\eqref{eqd}, whose constituents are
\begin{align*}
\tau(\phi) &= \tfrac{4\epsilon}{1+ \epsilon |x|^2} x ; \quad |\tau(\phi)|^2 = \tfrac{4.16}{(1+ \epsilon |x|^2)^4} |x|^2 ;\\
\Delta \lambda^2 &= \tfrac{-16\epsilon}{(1+ \epsilon |x|^2)^4} (4-2\epsilon|x|^2) ; \quad
\di \langle d\phi ,\tau(\phi)\rangle  = \tfrac{16\epsilon}{(1+ \epsilon |x|^2)^4} (4-2\epsilon|x|^2) ;\\
\nabla |\tau(\phi)|^2 &= \tfrac{8.16(1-3\epsilon|x|^2)}{(1+\epsilon|x|^2)^5} x ;
\quad \langle d\phi ,\tau(\phi)\rangle (\nabla |\tau(\phi)|^2) = \tfrac{8.16^2
(1-3\epsilon|x|^2)}{(1+\epsilon|x|^2)^8} |x|^2 ,
\end{align*}
therefore
\begin{align*}
&|\tau(\phi)|^{4} - 2\Delta \lambda^{2} |\tau(\phi)|^{2} + 4 \Delta \lambda^{2}
\di{\langle d\phi , \tau(\phi)\rangle} + 4 (\Delta \lambda^{2})^{2} +
2 \langle d\phi , \tau(\phi) \rangle (\nabla |\tau(\phi)|^{2}) \\
&= \tfrac{16^3}{(1+\epsilon|x|^2)^8} |x|^2 (1+2 \epsilon -3\epsilon |x|^2),
\end{align*}
so $\phi$ cannot be a biharmonic morphism.
\end{proof}

\begin{prop}
The identity map from $(\bn^n , ds^2)$ to $(\bn^n , h)$, where $\bn^n$ is the unit
ball in $\rn^n$, $ds^2$ its Euclidean metric and $h_{x} = \tfrac{4}{(1 - |x|^2)^2}
ds^2$ gives the hyperbolic space, is a biharmonic map if and only if $n=4$.
Furthermore, in either case, the biharmonic map is not a biharmonic morphism.
\end{prop}

\begin{proof}
Call $\phi$ the identity from $(\bn^n , ds^2)$ to $(\bn^n , h)$. Clearly $h =
e^{2\rho} ds^2$ for $\rho = \ln \tfrac{2}{1 - |x|^2}$ and
\begin{align*}
\nabla \rho &= \tfrac{2}{1 - |x|^2}x ; \quad
|\nabla \rho|^2 = \tfrac{4|x|^2}{(1 - |x|^2)^2} ;\\
\nabla|\nabla \rho|^2 &= 8\tfrac{1 + |x|^2}{(1 - |x|^2)^3}x ; \quad
\Delta \rho = \tfrac{2}{(1 - |x|^2)^2}( n(1 - |x|^2) +2 |x|^2) ;\\
\tr_{g} \nabla^{2} \nabla \rho &= \tfrac{4x}{(1 - |x|^2)^3} (n+2 - (n-2) |x|^2) .
\end{align*}
So Equation~\eqref{eqf} becomes
$$ \tfrac{8x}{(1 - |x|^2)^3} (4-n)(1+ |x|^2) =0 $$
and $\phi$ is biharmonic only in dimension four.\\
However, in dimension four, the Laplacian of its conformal factor $\lambda^2 =
\tfrac{4}{(1 - |x|^2)^2}$ is
$$ \Delta \lambda^2 = \tfrac{16}{(1 - |x|^2)^4} (4 + 2  |x|^2)$$
and
$$\tau(\phi) = -2 \g \rho = \tfrac{-4}{1 - |x|^2} x \, \mbox{ and } \quad
|\tau(\phi)|^2 = \tfrac{64|x|^2}{(1 - |x|^2)^4}  .$$
On the other hand:
$$ \langle d\phi , \tau(\phi) \rangle = \tfrac{-16}{(1 - |x|^2)^3} x \, \mbox{ and } \quad
\di \langle d\phi , \tau(\phi) \rangle = \tfrac{-16}{(1 - |x|^2)^4} (4 + 2 |x|^2) .$$
Therefore
$$ 4 \Delta \lambda^2 + 2 \di \langle d\phi , \tau(\phi) \rangle =
\tfrac{32}{(1 - |x|^2)^4} (4 + 2 |x|^2) \neq |\tau(\phi)|^2,$$ so $\phi$ is not a
biharmonic morphism.
\end{proof}

\begin{prop}
The identity map from $(\rn^{n}_{+}, ds^2)$ to $\hn^n = (\rn^{n}_{+},
\tilde{h}=\tfrac{1}{(x^{n})^2}ds^2)$ from the upper-half Euclidean space to the
hyperbolic space is biharmonic if and only if
$n=4$. Furthermore, in either case, the biharmonic map is not a biharmonic morphism.
\end{prop}

\begin{proof}
Call $\phi$ the identity map from $(\rn^{n}_{+}, ds^2)$ to $\hn^n$. Since $\tilde{h}
= e^{2\rho}$ for $\rho = -\ln x^{n}$ and
\begin{align*}
\g \rho &= -\tfrac{1}{x^{n}}\tfrac{\partial}{\partial x^{n}} ; \quad  |\g \rho|^2 = \tfrac{1}{(x^{n})^2} ;\\
\g |\g \rho|^2 &= -\tfrac{2}{(x^{n})^3}\tfrac{\partial}{\partial x^{n}} ; \quad \Delta
\rho = \tfrac{1}{(x^{n})^2},
\end{align*}
and Equation~\eqref{eqf} becomes
$$\tfrac{2}{(x^{n})^3}(n-4) \tfrac{\partial}{\partial x^{n}}=0.$$
On the other hand, in dimension four, its dilation is $\lambda^2 =
\tfrac{1}{(x^{4})^2}$ and $\tau(\phi) = -2 \g \rho =
\tfrac{2}{x^{4}}\tfrac{\partial}{\partial x^{4}}$, so
\begin{align*}
|\tau(\phi)|^2 &= \tfrac{4}{(x^{4})^4} ; \quad \Delta \lambda^2 = \tfrac{6}{(x^{4})^4} ;\\
\langle d\phi , \tau(\phi) \rangle &= \tfrac{2}{(x^{4})^3}\tfrac{\partial}{\partial
x^{4}} ; \quad \di \langle d\phi , \tau(\phi) \rangle = -\tfrac{6}{(x^{4})^4},
\end{align*}
so
$$ 4\Delta \lambda^2 + 2 \di \langle d\phi , \tau(\phi) \rangle = \tfrac{12}{(x^{4})^4} \neq |\tau(\phi)|^2$$
and $\phi$ is not a biharmonic morphism.
\end{proof}

Not all conformal maps in dimension four are biharmonic.
\begin{ex}
1) Consider the identity map $\phi$ from $(\hn^4 , g)=
(\rn_{+}^{4},\tfrac{1}{(x^4)^2}ds^2)$ to $(\rn^{4}_{+},ds^2)$. Then $ds^2 = e^{2\rho} g$
for $\rho = \ln x^4$ and, using the orthonormal basis $\{e_{i} = x^4
\tfrac{\partial}{\partial x^i}\}_{i=1,\dots,4}$, we have
\begin{align*}
\g^g \rho &= e_4 ; \quad |\g^g \rho|_{g}^2 = 1 ; \quad \g^g|\g^g \rho|^2 = 0 ;\\
\Delta^g \rho & = -3 ; \quad \tr_{g} \nabla^2 \g \rho = -3 e_{4}; \quad \ricci(\g^g \rho)  = -3 e_4,
\end{align*}
since $g= e^{2\alpha} ds^2$ for $\alpha = - \ln x^4$ and
\begin{align*}
&\nabla^{g}_{e_{1}} e_{4} = - e_{1} ; \quad \nabla^{g}_{e_{2}} e_{4} = - e_{2} ;\quad
\nabla^{g}_{e_{3}} e_{4} = - e_{3} ; \quad \nabla^{g}_{e_{4}} e_{4} = 0 ; \\
&\nabla^{g}_{e_{1}} e_{1} =  \nabla^{g}_{e_{2}} e_{2} = \nabla^{g}_{e_{3}} e_{3} =
e_{4}.
\end{align*}
Testing Equation~\eqref{eqf}, we have $-2 e_{4} \neq 0$, so $\phi$ is not biharmonic.\\
Moreover, from simple considerations
\begin{align*}
\tau(\phi) &= -2e_{4} ; \quad  \lambda^2 = (x^4)^2 ; \quad \Delta \lambda^2 = -2 (x^4)^2 ; \\
\langle d\phi , \tau(\phi) \rangle &= -2 e_{4} ; \quad \di \langle d\phi , \tau(\phi)
\rangle = 6 ,
\end{align*}
so $4\Delta\lambda^2 + 2 \di \langle d\phi , \tau(\phi) \rangle \neq |\tau(\phi)|^2$ and $\phi$ does not satisfy \eqref{eqd}.\\

2) Let $\phi$ be the identity map from $(\rn^n , g=\tfrac{4}{(1+\epsilon |x|^2)^2}
ds^2)$ to $(\rn^n , ds^2)$ ($\epsilon = \pm 1$). It is clearly conformal of dilation
$\lambda^2 = \tfrac{(1+\epsilon |x|^2)^2}{4}$ and $g = e^{2\rho} ds^2$ for $\rho = \ln 2 - \ln (1+\epsilon |x|^2)$. Then
$e_{i} = \tfrac{1+\epsilon |x|^2}{2}\tfrac{\partial}{\partial x^i}$ is an orthonormal basis for $g$ and
\begin{align*}
\nabla^{g}_{e_{i}}e_{j} &= \nabla^{ds^2}_{e_{i}}e_{j} + e_{i}(\rho)e_{j} + e_{j}(\rho)e_{i} - \langle
e_{i},e_{j} \rangle \g \rho  = - \epsilon x^j e_{i} + \delta_{ij} \epsilon X
\end{align*}
where $X= \sum_{k=1}^{n} x^k e_k$.
Therefore
\begin{align*}
&\g^{g} \ln\lambda = \tfrac{(1+\epsilon |x|^2)^2}{4}\tfrac{\partial}{\partial x^i}
(\ln\lambda) \tfrac{\partial}{\partial x^i} =  \epsilon X ;\\
&|\g^g \ln\lambda|_{g}^{2} = \tfrac{4}{(1+\epsilon |x|^2)^2} \tfrac{(1+\epsilon |x|^2)^2}{4} |x|^2 = |x|^2 ;\\
&\g^g |\g^g \ln\lambda|_{g}^{2} =\tfrac{(1+\epsilon |x|^2)^2}{4}
\tfrac{\partial}{\partial x^i} (|x|^2) \tfrac{\partial}{\partial x^i}
=(1+\epsilon |x|^2) X .
\end{align*}
Since $\Delta^{g} \ln\lambda = \sum_{i=1}^{4} e_{i}(e_{i} (\ln\lambda))
-(\nabla^{g}_{e_{i}}e_{i})(\ln\lambda)$ and
\begin{align*}
& \sum_{i=1}^{n} \nabla^{g}_{e_{i}}e_{i} = (n-1) \epsilon X ; \quad  X(\ln \lambda) =
\epsilon |x|^2 ; \quad
\sum_{i=1}^{n} (\nabla^{g}_{e_{i}}e_{i}) (\ln \lambda) = (n-1) |x|^2 ;\\
&e_{i} (\ln\lambda) = \epsilon x^i ;\quad  e_{i}(e_{i} (\ln\lambda)) = \epsilon
\tfrac{1+\epsilon |x|^2}{2} ,
\end{align*}
we have $\Delta^{g} \ln\lambda = \epsilon \tfrac{n + (2-n)\epsilon |x|^2}{2}$.\\
On the other hand
$$\Delta \g^g \ln\lambda = - \sum_{i=1}^{n} \nabla^{g}_{e_{i}}\nabla^{g}_{e_{i}} (\g^g \ln\lambda)
 - \nabla^{g}_{\nabla^{g}_{e_{i}}e_{i}} (\g^g \ln\lambda) , $$
and
\begin{align*}
& - \sum_{i=1}^{n} \nabla^{g}_{\nabla^{g}_{e_{i}}e_{i}} (\g^g \ln\lambda) = -
\nabla^{g}_{(n-1) \epsilon X} (\epsilon X)
= - (n-1) \tfrac{1+\epsilon |x|^2}{2} X ;\\
&\nabla^{g}_{e_{i}} (\g^g \ln\lambda)  = \epsilon (\tfrac{1-\epsilon |x|^2}{2} e_i + \epsilon x^i X) ;\\
&\nabla^{g}_{e_{i}}\nabla^{g}_{e_{i}} (\g^g \ln\lambda) = - \tfrac{1 + \epsilon |x|^2}{2} x^i e_i + (1+ \epsilon (x^i)^2)X ;\\
&\sum_{i=1}^{n} \nabla^{g}_{e_{i}}\nabla^{g}_{e_{i}} (\g^g \ln\lambda) = \tfrac{2n-1
+ \epsilon |x|^2}{2} X,
\end{align*}
so
$$- \Delta \g^g \ln\lambda = \tfrac{n + (2-n)\epsilon |x|^2}{2} X.$$
Equation~\eqref{eqf} becomes
$$ (2 + (4-n)\epsilon|x|^2) X =0,$$
which is impossible, whatever the value of $n$.
\end{ex}

\begin{ex}
Let $(M^2 ,h)$ be a Riemannian surface of Gaussian curvature $G_h$, $\beta : M^2 \times {\rn} \to {\rn}^{*}$ and
$\lambda : \rn \to \rn^{*}$ two positive functions. Consider the doubly twisted
product
$( {\rn} \times M^2 , g= \beta^2 dt^2 + \lambda^{-2} h)$.\\
Then, the projection
\begin{align*}
\phi  : ( {\rn} \times M^2 , g= \beta^2 dt^2 + \lambda^{-2} h) &\to (M^2 ,h) \\
(t,x) &\mapsto x
\end{align*}
is a biharmonic map if and only if
\begin{align}\label{2D}
0 &= -\lambda^4 \g_{h} (\Delta_{h} \ln\beta)  - 2\lambda^4 G_{h} \g_{h} \ln\beta +
\tfrac{3}{2} V(\lambda^2|\g_{h}\ln \beta|^2) V\\\notag &- \lambda^4|\g_{h}\ln
\beta|^2 \g_{h}\ln \beta - \nabla_{V}\nabla_{V} [\lambda^2 \g_{h}\ln \beta] -
\tfrac{\lambda^{4}}{2} \g_{h} |\g_{h} \ln\beta|_{h}^{2} \\\notag & - \lambda^2
V(V(\ln\lambda)\g_h \ln\beta - 2V(\ln\lambda)[\nabla_{V} (\lambda^2 \g_{h}\ln \beta)
- \lambda^2 |\g_h \ln\beta |_{h}^2 V].
\end{align}
In particular, for $\beta = c_2 \, e^{\int f(x) \, dx}$ with
$$f(x) = \frac{-c_1 (1+e^{c_1 x})}{1-e^{c_1 x}}$$
and $c_1 , c_2 \in \rn^*$, we have a family of Riemannian submersions
\begin{align*}
\phi : (\rn^2 \times \rn , dx^2 + dy^2 + \beta^2(x) dz^2) &\to (\rn^2 , dx^2 + dy^2)\\
\phi(x,y,z) &= (x,y)
\end{align*}
which are proper biharmonic maps.\\

In fact, it is easily checked that $\phi$ is a horizontally conformal
submersion of dilation $\lambda^2$.\\
Let  $V$ be the unit vertical vector $\frac{1}{\beta} \frac{d}{dt}$, using the Koszul
formula one can show that for vector fields $E$ and $F$  on ${\rn} \times M^2$
$$A_{E}F = g(E^\h,F^\h) V(\ln \lambda) V.$$
On the other hand, if $W$ is a vertical vector field, then
$$ B_{V} W = - \lambda^{2} \g_{h} (\ln\beta) \langle W,V\rangle .$$
Note also that the mean curvature of the fibres is $\mu = \nabla^{g}_{V} V = -
\lambda^{2} \g_{h} (\ln\beta)$ and the mean curvature of the horizontal distribution
is $\nu = V(\ln\lambda )
V$.\\

 Choosing a geodesic frame $\{ e_1 ,e_2 \}$ around a point
$p\in M^2$ and evaluating all subsequent formulas at a point $(p,t)\in M^2 \times
\rn$, straightforward computations show that (summing on repeated indices)
\begin{align*}
&[ \tr  ( (\nabla A)^{*}_{\mu} - (\nabla A)\mu + 3 A^{*}_{\nabla \mu} +
(\nabla_{\mu}B^{*})^{*} + 2 B_{\nabla_{\mu}} - B_{\nabla\mu} -2 A^{*}_{\nabla_{\mu}}) + 2 A^{*}_{\mu}\nu ]^\h \\
&= V(V(\ln\lambda))\mu + 3(V(\ln\lambda)^2)  \mu -3V(\ln\lambda) V(\lambda
e_{i}(\ln\beta)) \lambda e_{i}  \\
&+ \mu(g(\lambda e_{j},\mu)) \lambda e_{j} + |\mu|^{2} \mu .
\end{align*}

Since $\g^{\h} \ln\lambda =0$, the only remaining terms are
 \begin{align*}
[\nabla_{\mu}\mu]^\h &=  \tfrac{\lambda^4 }{2} \g_{h} ( | \g_{h} \ln\beta |_{h}^{2}) ;\\
[\nabla_{\mu} \nu]^\h &= - (V(\ln\lambda))^2 \mu ; \\
  \Delta^{\h} (\ln\lambda) &= -2(V(\ln\lambda))^2 ;\\
[\nabla_{\mu} \g\ln\lambda]^\h & = - (V(\ln\lambda))^2 \mu ;\\
 [\Delta \mu ]^\h &= \lambda^4 e_{i}(e_{i}(e_{j}\ln\beta))) e_{j} - (V(\ln\lambda))^2 \g^\h \ln\beta\\
&- V(\ln\lambda) V(\lambda e_{j}(\ln\beta)) \lambda e_{j} +  V(V(\lambda e_{j}(\ln\beta))) \lambda e_{j}\\
& - \lambda^4 |\g_{h} \ln\beta|_{h}^{2} \g_{h} \ln\beta 
+ \tfrac{\lambda^4}{2} \g_{h}(|\g_{h} \ln\beta|_{h}^{2}) \\
&+ \lambda e_{k}(\ln\beta)\lambda^2 [\nabla_{\lambda e_i}(\nabla^{M^2}_{e_i}e_k)]^\h .
\end{align*}
So the projection $\phi$ is biharmonic if and only if (still summing on repeated
indices)
\begin{align*}
0 &= -\lambda^4 e_{i}(e_{i}(e_{j}\ln\beta))) e_{j} - V(V(\lambda e_{j}(\ln\beta)))
\lambda e_{j} - \tfrac{\lambda^{4}}{2} \g_{h} |\g_{h} \ln\beta|_{h}^{2}\\
& - \lambda^4 G_{h} \g_{h} \ln\beta  
- \lambda^2 V(V(\ln\lambda)\g_h \ln\beta  - 2V(\ln\lambda)V(\lambda e_{i}(\ln\beta))
\lambda e_{i} \\
&- \lambda e_{k}(\ln\beta)\lambda^2 [\nabla_{\lambda
e_i}(\nabla^{M^2}_{e_i}e_k)]^\h .
\end{align*}
Using
\begin{align*}
&V(V(\lambda e_{j}(\ln\beta))) \lambda e_{j}  = - \tfrac{3}{2} V(\lambda^2
|\g_{h}\ln \beta|^2) V + \lambda^4|\g_{h}\ln \beta|^2 \g_{h}\ln \beta \\
& + \nabla_{V}\nabla_{V} [\lambda^2 \g_{h}\ln \beta] ;\\
&V(\lambda e_{j}(\ln\beta)) \lambda e_{j} = \nabla_{V} (\lambda^2 \g_{h}\ln \beta) - \lambda^2 |\g_h \ln\beta |_{h}^2 V ;\\
&\g_{h} (\Delta_{h} \ln\beta) = e_{i}e_{i}e_{j}(\ln\beta)e_{j} - G_h \g_{h}\ln\beta\\
&- \langle \nabla_{e_{1}}\nabla_{e_{1}}e_{1}+ \nabla_{e_{2}}\nabla_{e_{2}}e_{1}
,e_{2}\rangle e_{2}(\ln\beta)e_{1}
 - \langle \nabla_{e_{1}}\nabla_{e_{1}}e_{2} + \nabla_{e_{2}}\nabla_{e_{2}}e_{2}
,e_{1}\rangle e_{1}(\ln\beta) e_{2} ;\\
& \lambda e_{k}(\ln\beta)\lambda^2 [\nabla_{\lambda e_i}(\nabla^{M^2}_{e_i}e_k)]^\h
= -\lambda^4 [ e_{1}(\ln\beta) h(\nabla^{M^2}_{e_1}\nabla^{M^2}_{e_1}e_2
+ \nabla^{M^2}_{e_2}\nabla^{M^2}_{e_2}e_2,e_1) e_2 \\
&+ e_{2}(\ln\beta) h(\nabla^{M^2}_{e_1}\nabla^{M^2}_{e_1}e_1 +
\nabla^{M^2}_{e_2}\nabla^{M^2}_{e_2}e_1,e_2) e_1 ],
\end{align*}
we obtain the biharmonic equation (\ref{2D}) for $\phi$.\\

For the Riemannian submersion
\begin{align*}
\phi : (\rn^2 \times \rn , dx^2 + dy^2 + \beta^2(x) dz^2) &\to (\rn^2 , dx^2 + dy^2)\\
\phi(x,y,z) &= (x,y) ,
\end{align*}
Equation~\eqref{2D} reduces to
$$ ff' + f'' =0 ,$$
where $f = \ln \beta$. Solving this equation, we obtain the last statement in the
example.
\end{ex}

\end{document}